\documentclass{article}

	\newcommand{\mainTitle}{Weighted sums with two parameters of multiple zeta values and their formulas}
	
	\newcommand{\authorName}{MACHIDE, Tomoya}
	\newcommand{\organizationName}{Kinki University}
	\newcommand{\departmentName}{Interdisciplinary Graduate School of Science and Engineering}
	\newcommand{\majorName}{Research Center for Quantum Computing}
	\newcommand{\placeAddress}{3-4-1 Kowakae, Higashi-Osaka, Osaka 577-8502, Japan}
	\newcommand{\emailAddress}{E-mail: machide\_tomoya@alice.math.kindai.ac.jp, machide.t@gmail.com}
	\newcommand{\sectOne}{Introduction}
	\newcommand{\sectTwo}{Proof of Theorem \ref{1_Thm1}}
	\newcommand{\sectThree}{Proof of Theorem \ref{1_Thm2}}
	\newcommand{\sectFour}{Some identities derived from Theorem \ref{1_Thm1} and \ref{1_Thm2}}
\usepackage{geometry}                
\usepackage{ascmac}
\usepackage{amsmath}
\usepackage{amssymb}
\usepackage{amsthm}

	\newcommand{\nbk}[3]{#1#3#2}		
	\newcommand{\bgbk}[3]{\bigl{#1}#3\bigr{#2}}	
	\newcommand{\Bgbk}[3]{\Bigl{#1}#3\Bigr{#2}}			
	\newcommand{\bggbk}[3]{\biggl{#1}#3\biggr{#2}}			
	\newcommand{\Bggbk}[3]{\Biggl{#1}#3\Biggr{#2}}
	\newcommand{\autobk}[3]{\left#1#3\right#2}
	\newcommand{\mcbk}[4][?]{
						\ifx n#1
							\nbk{#2}{#3}{#4}
						\else
							\ifx b#1
								\bgbk{#2}{#3}{#4}
							\else
								\ifx B#1
									\Bgbk{#2}{#3}{#4}
								\else
									\ifx g#1
										\bggbk{#2}{#3}{#4}
									\else
										\ifx G#1
											\Bggbk{#2}{#3}{#4}
										\else
											\ifx a#1
												\autobk{#2}{#3}{#4}
											\else
												#4
											\fi
										\fi
									\fi
								\fi
							\fi
						\fi
					}
	\newcommand{\nsgsb}[1]{#1}		
	\newcommand{\bgsgsb}[1]{\big{#1}}	
	\newcommand{\Bgsgsb}[1]{\Big{#1}}			
	\newcommand{\bggsgsb}[1]{\bigg{#1}}			
	\newcommand{\Bggsgsb}[1]{\Bigg{#1}}
	\newcommand{\mcsgsb}[2][?]{
						\ifx n#1
							\nsgsb{#2}
						\else
							\ifx b#1
								\bgsgsb{#2}
							\else
								\ifx B#1
									\Bgsgsb{#2}
								\else
									\ifx g#1
										\bggsgsb{#2}
									\else
										\ifx G#1
											\Bggsgsb{#2}
										\else
											#2
										\fi
									\fi
								\fi
							\fi
						\fi
					}
	
	\newcommand{\setZ}[1][?]{\ifx?#1 \mathbb{Z} \else \mathbb{Z}_{#1} \fi}	
	\newcommand{\setQ}{\mathbb{Q}}	
		
	\newcommand{\setC}{\mathbb{C}}
	
	\newcommand{\gpSym}[2][n]{\ifx g#1\mathfrak{S}_{#2} \else S_{#2} \fi}
	\newcommand{\gpAlt}[2][n]{\ifx g#1 \mathfrak{A}_{#2} \else A_{#2} \fi }

	\newcommand{\vPack}[1][10]{\vspace{-#1pt}}
	
	\newcommand{\lnA}[1][]{&  &}
	\newcommand{\lnP}[2][n]{\ifx n#1\ \ #2 \ \ \else \ifx s#1 \ #2\ \else #2\fi\fi}
	\newcommand{\lnAP}[2][]{& #2 &}

	\newcommand{\lnAHs}[2][\nonumber]{#1 \\ & &\hspace{#2pt}}

	\newcommand{\refEq}[1]{(\ref{#1})}	
	\newcommand{\refSect}[1]{\S\ref{#1}}	
	\newcommand{\refH}[2]{#1 \ref{#2}}

	\newcommand{\bkR}[2][a]{\mcbk[#1]{(}{)}{#2}}
	\newcommand{\bkS}[2][a]{\mcbk[#1]{[}{]}{#2}}
	\newcommand{\bkB}[2][a]{\mcbk[#1]{\{}{\}}{#2}}

	\newcommand{\tmR}[4][?]{\bkR[#1]{#2} \ifx 0#4 \else \ifx 1#4 #3 \else {#3}^{#4} \fi \fi}	
	\newcommand{\tmS}[4][?]{\bkS[#1]{#2} \ifx 0#4 \else \ifx 1#4 #3 \else {#3}^{#4} \fi \fi}
	\newcommand{\tmM}[4][?]{\bkB[#1]{#2} \ifx 0#4 \else \ifx 1#4 #3 \else {#3}^{#4} \fi \fi}

	\newcommand{\undBT}[3][d]{ \underset{#3}{\underbrace{\if d#1 #2,\ldots,#2 \else #2 \fi}} }

	\newcommand{\Gam}{\Gamma}

	\newcommand{\mo}{(-1)}
	
	\newcommand{\opF}[3][?]{\ifx s#1 #2/#3 \else \ifx d#1 \dfrac{#2}{#3} \else \frac{#2}{#3} \fi\fi}

	\newcommand{\pw}[3][?]{\ifx 0#3 1  \else\ifx 1#3 \bkR[#1]{#2}\else \ifx?#1\bkR[#1]{#2}^{#3}\else {\bkR[#1]{#2}}^{#3} \fi\fi\fi}
	\newcommand{\id}[2]{#1_{#2}}
	\newcommand{\ip}[4][?]{\ifx ?#1 \ifx 1#4 {#2}_{#3} \else {#2}_{#3}^{#4} \fi \else \ifx n#1\ifx 1#3 {(#2_{#3})} \else (#2_{#3})^{#4} \fi \else #1\fi\fi}
	\newcommand{\pwR}[3][a]{ \pw[#1]{#2}{#3} }
	\newcommand{\pwB}[3][a]{\ifx 0#3 1  \else\ifx 1#3 \bkB[#1]{#2}\else \ifx?#1\bkB[#1]{#2}^{#3}\else {\bkB[#1]{#2}}^{#3} \fi\fi\fi}
					
	\newcommand{\nFc}[3][n]{#2\bkR[#1]{#3}}
	
	\newcommand{\idFc}[4][n]{\id{#2}{#3}\bkR[#1]{#4}}
	\newcommand{\pwFc}[4][n]{\pw{#2}{#3}\bkR[#1]{#4}}
	\newcommand{\ipFc}[5][n]{\ip{#2}{#3}{#4}\bkR[#1]{#5}}
		\newcommand{\Fc}{\nFc}
	
	\newcommand{\ptl}{\partial}

	\newcommand{\pDer}[2][?]{\opF[#1]{\ptl}{\ptl#2}}	
	\newcommand{\bpDer}[2][?]{\bkR{\opF[#1]{\ptl}{\ptl#2}}}	
	\newcommand{\pDerT}[3][?]{\opF[#1]{\ptl^{#3}}{#2}}

	\newcommand{\tpT}[3][a]{ {#2}\atop \bkR[#1]{#3} }

	\newcommand{\nSm}[2][?]{\ifx l#1 \sum\limits_{#2} \else \sum_{#2} \fi}
	\newcommand{\nSmT}[3][?]{\ifx l#1 \sum\limits_{#2}^{#3} \else \sum_{#2}^{#3} \fi}	
	\newcommand{\pSm}[2][?]{\ifx t#1 \sum_{#2}^{\prime} \else \sideset{}{^\prime}\sum_{#2} \fi}
	\newcommand{\pSmT}[3][?]{\ifx t#1 \sum_{#2}^{\prime#3} \else \sideset{}{^\prime}\sum_{#2}^{#3} \fi}	
	\newcommand{\tpTSm}[3][?]{ \nSm[#1]{ \tpT{#2}{#3} } }
	\newcommand{\tpTSmT}[4][?]{ \nSmT[#1]{ \tpT{#2}{#3} }{#4} }
	
		\newcommand{\Sm}{\nSm}
		\newcommand{\SmT}{\nSmT}
		\newcommand{\tpSm}{\tpTSm}
		\newcommand{\tpSmT}{\tpTSmT}
	\newcommand{\nPd}[2][?]{\ifx l#1 \prod\limits_{#2} \else \prod_{#2} \fi}
	\newcommand{\nPdT}[3][?]{\ifx l#1 \prod\limits_{#2}^{#3} \else \prod_{#2}^{#3} \fi}	
	
	\newcommand{\tpTPdT}[4][?]{ \nPdT[#1]{ \tpT{#2}{#3} }{#4} }
		
		\newcommand{\PdT}{\nPdT}
		
		\newcommand{\tpPdT}{\tpTPdT}

	\newcommand{\nLim}[2][?]{\ifx l#1 \lim\limits_{#2} \else \lim_{#2} \fi}
		\newcommand{\Lim}{\nLim}

	\newcommand{\glcondEnvLineHead}[1]{ \ifx*#1 \begin{eqnarray*} \else \begin{eqnarray}  \label{#1} \fi }
	\newcommand{\glcondEnvLineTail}[1]{ \ifx*#1 \end{eqnarray*} \else \end{eqnarray} \fi }
	\newcommand{\glcondDis}[1]{\ifx d#1 \displaystyle \fi}

		\newcommand{\envHLineT}[3][*]{ \glcondEnvLineHead{#1} #2&=&#3\glcondEnvLineTail{#1} }
		\newcommand{\envHLineTDef}[3][*]{ \glcondEnvLineHead{#1} #2&:=&#3\glcondEnvLineTail{#1} }
			\newcommand{\envHLine}{\envHLineT}
			\newcommand{\envHLineDef}{\envHLineTDef}

		\newcommand{\envPLine}[2][*]{\glcondEnvLineHead{#1} #2\glcondEnvLineTail{#1}}
		
		\newcommand{\envMO}[2][*]{$\ifx d#1 \displaystyle \fi#2$}
		\newcommand{\envMT}[3][*]{$\ifx d#1 \displaystyle \fi#2=#3$}
		\newcommand{\envMTDef}[3][*]{$\ifx d#1 \displaystyle \fi#2:=#3$}
			\newcommand{\envM}{\envMT}
			
		\newcommand{\envMTh}[4][*]{$\ifx d#1 \displaystyle \fi#2=#3=#4$}
		\newcommand{\envMF}[5][*]{$\ifx d#1 \displaystyle \fi#2=#3=#4=#5$}
		
	\newcommand{\envMLineT}[3][*]{ \ifx*#1 \begin{multline*} #2\lnP{=}#3\end{multline*} \else \begin{multline} \label{#1} #2\lnP{=}#3\end{multline} \fi }

	\newcommand{\abs}[1]{\left | #1 \right |  }		
	
	\newcommand{\sbLandau}[2][a]{\Fc[#1]{O}{#2}}
	
	\newcommand{\sbPochD}[3][a]{\bkR[#1]{#2}_{#3}}
		
	\newcommand{\exx}[2][B]{ \Fc[#1]{\exp}{#2} }

	\newcommand{\tfcS}[2][?]{\Fc[#1]{\sin}{#2}}
	
	\newcommand{\tfcT}[2][?]{\Fc[#1]{\tan}{#2}}
	\newcommand{\tfcCT}[2][?]{\Fc[#1]{\cot}{#2}}

	\newcommand{\lgg}[2][?]{\Fc[#1]{\log}{#2}}

	\newcommand{\fcZeta}[1]{\zeta(#1)}
		\newcommand{\fcZ}{\fcZeta}		
	\newcommand{\lgP}[3][n]{\idFc[#1]{Li}{#2}{#3}}
	
	\newcommand{\fcGam}[2][n]{\Fc[#1]{\Gam}{#2}}

	\newcommand{\glcondEnvLineTailPd}[1]{.\ifx*#1 \end{eqnarray*} \else \end{eqnarray} \fi }
	\newcommand{\glcondEnvLineTailCm}[1]{,\ifx*#1 \end{eqnarray*} \else \end{eqnarray} \fi }

	\newcommand{\envProof}[2][?]{ \par\vspace{-5pt}\mbox{}\\ \ifx?#1\emph{Proof.}\else\emph{Proof of #1.}\fi \ #2 \hfill $\Box$\\ \par}
	
		\newcommand{\envLineTPd}[3][*]{ \glcondEnvLineHead{#1} & &#2\\&=&#3\nonumber \glcondEnvLineTailPd{#1} }
		
		\newcommand{\envLineTCm}[3][*]{ \glcondEnvLineHead{#1} & &#2\\&=&#3\nonumber \glcondEnvLineTailCm{#1} }
		
			\newcommand{\envLinePd}{\envLineTPd}
			
			\newcommand{\envLineCm}{\envLineTCm}
			
		\newcommand{\envLineThPd}[4][*]{ \glcondEnvLineHead{#1} & &#2\\&=&#3\nonumber \\&=&#4\nonumber \glcondEnvLineTailPd{#1} }
		\newcommand{\envLineThCm}[4][*]{ \glcondEnvLineHead{#1} & &#2\\&=&#3\nonumber \\&=&#4\nonumber \glcondEnvLineTailCm{#1} }
		
		\newcommand{\envLineFCm}[5][*]{ \glcondEnvLineHead{#1} & &#2\\&=&#3\nonumber \\&=&#4\nonumber \\&=&#5\nonumber \glcondEnvLineTailCm{#1} }
		
		\newcommand{\envLineFiCm}[6][*]{ \glcondEnvLineHead{#1} & &#2\\&=&#3\nonumber\\&=&#4\nonumber\\&=&#5\nonumber\\&=&#6\nonumber\glcondEnvLineTailCm{#1} }

		\newcommand{\envLineSePd}[8][*]{ \glcondEnvLineHead{#1} & &#2\\&=&#3\nonumber\\&=&#4\nonumber\\&=&#5\nonumber\\&=&#6\nonumber\\&=&#7\nonumber\\&=&#8\nonumber\glcondEnvLineTailPd{#1}	 }
		
		\newcommand{\envHLineTPd}[3][*]{ \glcondEnvLineHead{#1} #2&=&#3\glcondEnvLineTailPd{#1} }
		\newcommand{\envHLineTDefPd}[3][*]{ \glcondEnvLineHead{#1} #2&:=&#3\glcondEnvLineTailPd{#1} }
		\newcommand{\envHLineTCm}[3][*]{ \glcondEnvLineHead{#1} #2&=&#3\glcondEnvLineTailCm{#1} }
		\newcommand{\envHLineTCmDef}[3][*]{ \glcondEnvLineHead{#1} #2&:=&#3\glcondEnvLineTailCm{#1} }
			\newcommand{\envHLinePd}{\envHLineTPd}
			\newcommand{\envHLineDefPd}{\envHLineTDefPd}
			\newcommand{\envHLineCm}{\envHLineTCm}
			\newcommand{\envHLineCmDef}{\envHLineTCmDef}

		\newcommand{\envHLineThCm}[4][*]{ \glcondEnvLineHead{#1} #2&=&#3\\&=&#4\nonumber\glcondEnvLineTailCm{#1}}
		
		\newcommand{\envHLineFCm}[5][*]{ \glcondEnvLineHead{#1} #2&=&#3\\&=&#4\nonumber \\&=&#5\nonumber \glcondEnvLineTailCm{#1} }

		\newcommand{\envHLineCFPd}[5][*]{\glcondEnvLineHead{#1} #2&=&#3,\\#4&=&#5\nonumber\glcondEnvLineTailPd{#1}}

		\newcommand{\envHLineCEPd}[9][*]{\glcondEnvLineHead{#1} #2&=&#3,\\#4&=&#5,\nonumber\\#6&=&#7,\nonumber\\#8&=&#9\nonumber\glcondEnvLineTailPd{#1}}

		\newcommand{\envPLineCm}[2][*]{\glcondEnvLineHead{#1} #2\glcondEnvLineTailCm{#1}}

		\newcommand{\envOTLinePd}[4][*]{\glcondEnvLineHead{#1} #2\lnAP{=}#3\lnP{=}#4.\glcondEnvLineTail{#1}}
		\newcommand{\envOTLineCm}[4][*]{\glcondEnvLineHead{#1} #2\lnAP{=}#3\lnP{=}#4,\glcondEnvLineTail{#1}}
		\newcommand{\envOTLineDefPd}[4][*]{\glcondEnvLineHead{#1} #2\lnAP{:=}#3\lnP{=}#4.\glcondEnvLineTail{#1}}

		\newcommand{\envMOPd}[2][*]{$\ifx d#1 \displaystyle \fi#2.$}
		\newcommand{\envMOCm}[2][*]{$\ifx d#1 \displaystyle \fi#2,$}
		\newcommand{\envMTPd}[3][*]{$\ifx d#1 \displaystyle \fi#2=#3.$}
		\newcommand{\envMTCm}[3][*]{$\ifx d#1 \displaystyle \fi#2=#3,$}
		\newcommand{\envMTDefPd}[3][*]{$\ifx d#1 \displaystyle \fi#2:=#3.$}
		\newcommand{\envMTDefCm}[3][*]{$\ifx d#1 \displaystyle \fi#2:=#3,$}
			
			\newcommand{\envMCm}{\envMTCm}

	\newcommand{\envELineT}[3][*]{ \ifx*#1 \begin{equation*} #2\lnP{=}#3\end{equation*} \else \begin{equation} \label{#1} #2\lnP{=}#3\end{equation} \fi }
	\newcommand{\envELineTPd}[3][*]{ \ifx*#1 \begin{equation*} #2\lnP{=}#3.\end{equation*} \else \begin{equation} \label{#1} #2\lnP{=}#3.\end{equation} \fi }
	\newcommand{\envELineTCm}[3][*]{ \ifx*#1 \begin{equation*} #2\lnP{=}#3,\end{equation*} \else \begin{equation} \label{#1} #2\lnP{=}#3,\end{equation} \fi }

	\newcommand{\envMLineTPd}[3][*]{ \ifx*#1 \begin{multline*} #2\lnP{=}#3.\end{multline*} \else \begin{multline} \label{#1} #2\lnP{=}#3.\end{multline} \fi }
	\newcommand{\envMLineTCm}[3][*]{ \ifx*#1 \begin{multline*} #2\lnP{=}#3,\end{multline*} \else \begin{multline} \label{#1} #2\lnP{=}#3,\end{multline} \fi }
		\newcommand{\envMLinePd}{\envMLineTPd}
		\newcommand{\envMLineCm}{\envMLineTCm}

		\newcommand{\envHLineCEPdPt}[9]{\begin{eqnarray*}#2&#1&#3,\\#4&#1&#5,\nonumber\\#6&#1&#7,\nonumber\\#8&#1&#9.\nonumber\end{eqnarray*}}

			\newcommand{\OTLineCThCm}[4][?]{#2&=&#3=#4,\nonumber \ifx#1p \\ \fi}
			\newcommand{\OTLineCThPd}[4][?]{#2&=&#3=#4.\nonumber \ifx#1p \\ \fi}
			\newcommand{\OTLineCThCmDef}[4][?]{#2&:=&#3=#4,\nonumber \ifx#1p \\ \fi}
			\newcommand{\OTLineCThDefPd}[4][?]{#2&:=&#3=#4.\nonumber \ifx#1p \\ \fi}
			\newcommand{\OTLineCSCm}[7][?]{#2&=&#3=#4\nonumber#5&=&#6=#7,\nonumber \ifx#1p \\ \fi}
			\newcommand{\OTLineCSPd}[7][?]{#2&=&#3=#4\nonumber#5&=&#6=#7.\nonumber \ifx#1p \\ \fi}

	\newcommand{\envMLineTCmPt}[4][*]{ \ifx*#1 \begin{multline*} #3\lnP{#2}#4,\end{multline*} \else \begin{multline} \label{#1} #3\lnP{#2}#4,\end{multline} \fi }
	\newcommand{\envMLineTPdPt}[4][*]{ \ifx*#1 \begin{multline*} #3\lnP{#2}#4.\end{multline*} \else \begin{multline} \label{#1} #3\lnP{#2}#4.\end{multline} \fi }

\theoremstyle{plain}
\newtheorem{theorem}{THEOREM}[section]
\newtheorem{proposition}[theorem]{PROPOSITION}
\newtheorem{lemma}[theorem]{LEMMA}

\theoremstyle{definition}

\theoremstyle{remark}
\newtheorem{remark}[theorem]{REMARK}
\allowdisplaybreaks[3]
\numberwithin{equation}{section}

	\newcommand{\lccondBibitem}[3][]{ \if ?#2 \bibitem{#3} \else \bibitem[#2]{#3} \fi}
	\newcommand{\refPaper}[8][?]{
			\lccondBibitem{#1}{#2}
				#3,			
				\emph{#4}, 	
				#5\ 			
				{\bf #6},		
				#7,			
				#8.			
		}
	\newcommand{\refPreprint}[6][?]{
			\lccondBibitem{#1}{#2}
				#3,			
				\emph{#4}, 	
				preprint; #5,			
				#6.			
		}
	\newcommand{\refPaperRep}[9][?]{
			\lccondBibitem{#1}{#2}
				#3,			
				\emph{#4}, 	
				#5\ 			
				{\bf #6},		
				#7,			
				#8			
				; reprinted in #9	
		}

	\newcommand{\refPaperAlm}[5][?]{
			\lccondBibitem{#1}{#2}
				#3,	 		
				\emph{#4}, 	
				#5		
		}

	\newcommand{\refRem}[1]{\refH{Remark}{#1}}
\DeclareFontEncoding{OT2}{}{}
\DeclareFontSubstitution{OT2}{cmr}{m}{n}

\DeclareFontFamily{OT2}{cmr}{\hyphenchar\font45 }
\DeclareFontShape{OT2}{cmr}{m}{n}{%
   <5><6><7><8><9>gen*wncyr%
   <10><10.95><12><14.4><17.28><20.74><24.88>wncyr10}{}
\DeclareFontShape{OT2}{cmr}{b}{n}{%
   <5><6><7><8><9>gen*wncyb%
   <10><10.95><12><14.4><17.28><20.74><24.88>wncyb10}{}

\DeclareMathAlphabet{\mathcyr}{OT2}{cmr}{m}{n}
\DeclareMathAlphabet{\mathcyb}{OT2}{cmr}{b}{n}
\SetMathAlphabet{\mathcyr}{bold}{OT2}{cmr}{b}{n}

	\newcommand{\admInd}[3][-]{
					\ifx-#1	\Omega({#2};{#3})
					\else		\Omega_{#1}({#2};{#3})
					\fi			
				}

	\newcommand{\algFZV}[1][?]{\ifx ?#1 \mathfrak{H} \else \mathfrak{H}^{#1} \fi}
	
	\newcommand{\sumZV}[4][n]{{\mathfrak{Z}}_{#2}^{(#3)}\mcbk[#1]{(}{)}{#4}}

		\newcommand{\smZV}{\sumZV}

	\newcommand{\lsmDZV}[3][n]{\idFc[#1]{\mathfrak{Z}}{#2}{#3}}
	\newcommand{\lsmZV}[4][n]{\ipFc[#1]{\mathfrak{Z}}{#2}{(#3)}{#4}}
	\newcommand{\lsmZVd}[5][n]{\ipFc[#1]{\mathfrak{Z}}{#2}{(#3);(#4)}{#5}}
	\newcommand{\lcoeSmZV}[4][n]{\ipFc{\mathfrak{z}}{#2}{(#3)}{#4}}
	\newcommand{\lsmMP}[5][n]{\ipFc[#1]{\mathfrak{L}}{#2}{(#3)}{#4;#5}}
	
	\newcommand{\lgfcZV}[2][n]{\Fc[#1]{\mathfrak{Z}}{#2}}
	\newcommand{\lfcP}[3][n]{ \pwFc[#1]{P}{(#2)}{#3} } 
	\newcommand{\lfcF}[3][n]{ \idFc[#1]{f}{#2}{#3} } 

\allowdisplaybreaks[4]
	\title{\mainTitle}
	\author{\authorName}
	\date{}

\begin{document}
\maketitle
\begin{abstract}
A typical formula of multiple zeta values is the \emph{sum formula} which expresses a Riemann zeta value as 
	a sum of all multiple zeta values of fixed weight and depth.
Recently \emph{weighted sum formulas}, which are weighted analogues of the sum formula, have been studied by many people.
In this paper, we give two formulas of weighted sums with two parameters
	of multiple zeta values.
As applications of the formulas, 
	we find some linear combinations of multiple zeta values 
	which can be expressed as polynomials of usual zeta values 
	with coefficients in the rational polynomial ring generated by the two parameters,
	and obtain some identities for weighted sums of multiple zeta values of small depths.
\end{abstract}
\section{\sectOne}
Let $n$ be a positive integer.
For any multi-index $(l_1,l_2,\ldots,l_n)$ of positive integers with $l_1\geq2$,
	the \emph{multiple zeta value}, 
	sometimes called \emph{multiple harmonic series} or \emph{Euler sum},
	is a real number defined by the convergent series
	\envHLineDefPd[1_Def_MZV]
	{
		\fcZ{l_1,l_2,\ldots,l_n}
	}
	{
		\Sm{m_1>m_2>\cdots>m_n>0} \opF{1}{ \pw{m_1}{l_1}\pw{m_2}{l_2}\cdots\pw{m_n}{l_n} }
	}
We call the integer $l:=l_1+\cdots+l_n$ its weight and $n$ its depth.
These real numbers, already considered by Euler \cite{Euler1775} in the case of depth $2$, 
	have arisen in various areas of arithmetical algebraic geometry, number theory, knot theory and physics since the early 1990 (see \cite{Zagier94}).
A principal goal of the study of multiple zeta values is to determine the algebraic relations among them as much as possible \cite{Hoffman97,IKZ06}.
Recently Brown \cite{Brown11} has given a crucial result, or proved a conjecture proposed by Hoffman \cite{Hoffman97},
	which states that
	every multiple zeta value of weight $l$ is a $\setQ$-linear combination of $\fcZ{l_1,\ldots,l_n}$'s
	satisfying $n\geq1$, $l_1,\ldots,l_n\in\bkB{2,3}$ and $l=l_1+\cdots+l_n$.
The result yields the upper bound for the dimension of the $\setQ$-vector space 
	spanned by the multiple zeta values of any fixed weight $l$,
	which was conjectured by Zagier \cite{Zagier94} and proved independently 
	in the papers \cite{DG05} of Deligne and Goncharov and \cite{Terasoma02} of Terasoma.

A way of researching the algebraic relations among multiple zeta values is to give their elegant identities
	which were studied in many papers \cite{BBB97,ELO09,Hoffman92,HO03,Kawashima09,Ohno99,Tanaka09}.
Elegant identities are also interesting in themselves.
One of the typical examples of them is the \emph{sum formula} 
	\envHLine[1_Eq_SF]
	{
		\tpSm{l_1\geq2,l_2,\ldots,l_n\geq1}{l_1+\cdots+l_n=l} \fcZ{l_1,\ldots,l_n}
	}
	{
		\fcZ{l}
	}
	for any pair $(l,n)$ of positive integers with $l>n\geq1$, which means that the sum of all the multiple zeta values of fixed weight $l$ and depth $n$ 
	is equal to the special value of the Riemann zeta function at argument $l$.
This formula was conjectured in \cite{Hoffman92} and proved by Granville \cite{Granville97} and Zagier \cite{Zagier95}.
When $n=2$, it is the \emph{Euler sum formula} \cite{Euler1775}
	\envHLinePd[1_Eq_eulerSF]
	{
		\SmT{j=2}{l-1} \fcZ{j,l-j} 
	}
	{
		\fcZ{l}
	}
	
Recently weighted analogues of the sum formula have been discovered or studied \cite{Nakamura09,Sasaki11,SC12}.
Ohno and Zudilin \cite{OZu08} 
	gave the \emph{weighted Euler sum formula}
	\envHLine[1_Eq_eulerWSF]
	{
		\SmT{j=2}{l-1} \pw{2}{j}\fcZ{j,l-j} 
	}
	{
		(l+1)\fcZ{l}
	}
	to prove relations between multiple zeta values and \emph{multiple zeta star values}
	whose definition is same as \refEq{1_Def_MZV} except allowing equality among the $m_i$'s.
Gangl, Kaneko and Zagier \cite{GKZ06} systematically studied the double zeta values in a formal setting, and obtained many formulas.
One \cite[(26)]{GKZ06} of the formulas includes the weighted sum with parameters $x,y$ of the double zeta values of weight $l$
	\envHLineCmDef[1_Def_WSFp]
	{
		\lsmDZV{l}{x,y} 
	}
	{
		\SmT{j=2}{l-1} x^{j-1}y^{l-j-1} \fcZ{j, l-j}
	}
	which is the generating function of the double zeta values of weight $l$ in other words.
The formula we call the \emph{weighted sum formula with parameters} is
	\envHLinePd[1_Eq_WSFp]
	{
		\lsmDZV{l}{x+y,y} + \lsmDZV{l}{y+x,x} 
	}
	{
		\lsmDZV{l}{x,y} + \lsmDZV{l}{y,x} + \opF{x^{l-1}-y^{l-1}}{x-y}\fcZ{l}
	}
We note that this with $(x,y)=(1,0)$ and $x=y=1$ yield \refEq{1_Eq_eulerSF} and \refEq{1_Eq_eulerWSF} respectively,
	and that the coefficients of $x^{j-1}y^{l-j-1}$ in the formula are equivalent to the \emph{double shuffle relations} 
	derived from \emph{shuffle relations} and \emph{harmonic relations} among double zeta values (see \cite{GKZ06,IKZ06} for the relations).
Motivated by the work of Ohno and Zudilin,
	Guo and Xie \cite{GX09} generalized the weighted Euler sum formula \refEq{1_Eq_eulerWSF} to arbitrary depth $n$.
Eie, Yang and Ong \cite[Proposition 8]{EYO10} also obtained a formula
	which gives relations between weighted sums of multiple zeta values of even depth
	and sums of products of two multiple zeta values $\fcZ{k+1,\pwB{1}{j-1}}\ (k,j\geq1)$
	in order to prove alternating double sum formulas of multiple zeta values \cite[Theorem 4 and 5]{EYO10}.
Here $\pwB{1}{k}$ denotes $k$-tuple $(1,\ldots,1)$.

In this paper, we give two formulas including the
	following weighted sums with two parameters of multiple zeta values of any weight $l$ and any depth $n$,
	\envHLineCmDef[1_Def_wsmMZV]
	{
		\lsmZV{l}{n}{x,y}
	}
	{
		\tpSm{l_1\geq2,l_2,\ldots,l_n\geq1}{l_1+\cdots+l_n=l} \pw{x}{l_1-1}\pw{y}{l-l_1-(n-1)} \fcZ{l_1,\ldots,l_n}
	}
	where $\lsmZV{l}{1}{x,y}$ stands for $\pw{x}{l-1}\fcZ{l}$ and $\lsmZV{l}{2}{x,y}$ equals $\lsmDZV{l}{x,y}$.
One of the formulas is given by a generating function expression for simplicity.
\begin{theorem}\label{1_Thm1}
Let $\fcGam{z}$ be the gamma function.
For complex numbers $x$ and $y$, we have
	\envLineThPd[1_Thm1_Eq1]
	{
		\Sm{l>n\geq1} \bkS[g]{ \lsmZV{l}{n}{x+y,y} + \mo^{n} \lsmZV{l}{n}{y+x,x} - \bkR[b]{ \mo^{n}\pw{x}{l-n} + \pw{y}{l-n} } \fcZ{l} } \pw{X}{l-n} \pw{Y}{n}
	}
	{
		1 - \opF{ \fcGam{1-xX}\fcGam{1-Y} \fcGam{1-yX}\fcGam{1+Y} }{ \fcGam{1-xX-Y}\fcGam{1-yX+Y} }
	}
	{
		1 - \exx[a]{ \SmT{m=2}{\infty} \fcZ{m} \opF{ \bkR{x^m+y^m}X^m + \bkR{1+\mo^m}Y^m - \pwR{xX+Y}{m} - \pwR{yX-Y}{m} }{m} }
	}
\end{theorem}
Another formula gives a relation among 
	$\lsmZV{l}{n}{x+y,y} + \mo^{n} \lsmZV{l}{n}{y+x,x}\ (2\leq n\leq l-1)$ for any weight $l$, 
	which yields the weighted Euler sum formula \refEq{1_Eq_eulerWSF} of Ohno and Zudilin by substituting $(1,1)$ for $(x,y)$.
\begin{theorem}\label{1_Thm2}
For any integer $l\geq3$ and complex numbers $x,y$, we have
	\envLinePd[1_Thm2_Eq1]
	{
		\SmT{n=2}{l-1} \pwR{y-x}{n-2} \bkR{ \lsmZV{l}{n}{x+y,y} + \mo^{n} \lsmZV{l}{n}{y+x,x} }
	}
	{
		\bkS{ 
			\opF{\pw{y}{l+1}-\pw{x}{l+1}}{yx(y-x)}
			- 
			\pwR{y-x}{l-2} \bkR[g]{ (1+\mo^l) \opF{2^{l-1}-1}{2^{l-1}} + \mo^{l}\opF{x}{y} + \opF{y}{x} }
		}
		\fcZ{l}
	}
\end{theorem}
We comment about an application of \refH{Theorem}{1_Thm1}.
Many people 
	\cite{KO10,Li11,OZa01,Yamazaki10}
	recently have given some combinations among multiple zeta or zeta star values
	which can be written in terms of polynomials of single zeta values with rational coefficients.
We see from \refH{Theorem}{1_Thm1} that
	every $\lsmZV{l}{n}{x+y,y} + \mo^{n} \lsmZV{l}{n}{y+x,x}$ can be expressed as such a combination,
	but the coefficients of the polynomial have not only rational numbers but also the parameters $x$ and $y$.
	
We also give a generalization of the formula \cite[Proposition 8]{EYO10} of Eie, Yang and Ong
	in the course of the proof of Theorem \ref{1_Thm1},
	and show that the generalization is equivalent to the formula of Arakawa and Kaneko \cite[Corollary 11]{AK99}
	(see Proposition \ref{2_Prop1} and Remark \ref{2_Rem2}).
We note that
	\refH{Theorem}{1_Thm1} and \ref{1_Thm2} do not include the formula \refEq{1_Eq_WSFp} of Gangl, Kaneko and Zagier
	but \refH{Theorem}{1_Thm1} implies it by use of some equations of double zeta values (see \refRem{4_Rem1}),
	which suggests that each of \refEq{1_Thm1_Eq1} and \refEq{1_Thm2_Eq1} is a formula of a different type from \refEq{1_Eq_WSFp}.
	
The paper is organized as follows;
\refSect{sectTwo} and \refSect{sectThree} devote the proofs of Theorem \ref{1_Thm1} and \ref{1_Thm2} respectively.
In \refSect{sectFour}, we give some identities among multiple zeta values of small depths between $2$ and $4$ 
	as applications of the theorems.

\section{\sectTwo} \label{sectTwo}
For a proof of Theorem \ref{1_Thm1},
	we first consider weighted sums with two parameters of \emph{multiple polylogarithms} instead of multiple zeta values.
The multiple polylogarithm is defined by the convergent series
	\envHLineDef
	{
		\lgP{l_1,\ldots,l_n}{z}
	}
	{
		\Sm{m_1>\cdots>m_n>0} \opF{ \pw{z}{m_1} }{ \pw{m_1}{l_1}\cdots\pw{m_n}{l_n} }
	}
	for any multi-index $(l_1,\ldots,l_n)$ of positive integers and any complex number $z$ with $\abs{z}<1$.
The weighted sum with parameters of multiple polylogarithms of fixed weight $l$ is
	\envHLineCmDef
	{
		\lsmMP{l}{n}{x,y}{z}
	}
	{
		\tpSm{l_1,\ldots,l_n\geq1}{l_1+\cdots+l_n=l} \pw{x}{l_1-1}\pw{y}{l-l_1-(n-1)} \lgP{l_1,\ldots,l_n}{z}
	}
	where $l\geq n\geq1$ and $x,y,z\in\setC$ with $\abs{z}<1$.
The definition of $\lsmMP{l}{n}{x,y}{z}$ allows $l_1=1$ in the summand unlike the case of the multiple zeta values,
	because the multiple polylogarithms $\lgP{l_1,\ldots,l_n}{z}$ with $l_1=1$ are convergent by $\abs{z}<1$
	and needed in the following identities.
\begin{lemma}\label{2_Lem1}
For integers $l,n$ with $l\geq n\geq1$ and complex numbers $x,y,z$ with $\abs{z}<1$, we have
	\envLinePd[2_Lem1_Eq1]
	{
		\lsmMP{l}{n}{x+y,y}{z} +\mo^{n}\lsmMP{l}{n}{y+x,x}{z}
	}
	{
		\tpSm{j_1,j_2\geq1}{j_1+j_2=n} \mo^{j_2-1} \tpSm{k_1,k_2\geq0}{k_1+k_2=l-n}
		\pw{x}{k_1} \pw{y}{k_2} \lgP{k_1+1,\pwB{1}{j_1-1}}{z} \lgP{k_2+1,\pwB{1}{j_2-1}}{z}
	}
\end{lemma}
\envProof{
Let $m_1,\ldots,m_n$ be complex numbers.
We know the identity
	\envHLinePd
	{
		\PdT{i=1}{n} \opF{1}{m_i}
	}
	{
		\mo^{n-1} \SmT{j=1}{n} \opF{1}{m_j} \tpPdT{i=1}{i\neq j}{n} \opF{1}{m_j-m_i}
	}
For example, it was used in the paper \cite{Granville97} for the proof of the sum formula \refEq{1_Eq_SF}.
By putting
	\envHLineDef
	{
		\lfcP{n}{ \id{m}{1},\ldots,\id{m}{n} }
	}
	{
		\PdT{i=1}{n} \opF{1}{ m_i+\cdots+m_n }
	}
and replacing $m_i$ with $m_i+\cdots+m_n$ in the identity for $i=1,\ldots,n$, we get
	\envLineThPd
	{
		\lfcP{n}{ m_1,\ldots,m_n }
	}
	{
		\mo^{n-1} \SmT{j=1}{n} \opF{1}{m_j+\cdots+m_n} 
		\bkS[G]{ \PdT{i=1}{j-1} \opF{1}{-m_{j-1}-\cdots-m_i} } \bkS[G]{ \PdT{i=j+1}{n} \opF{1}{m_j+\cdots+m_{i-1}} }
	}
	{
		\SmT{j=1}{n} \mo^{n-j} \lfcP{ j-1 }{ \id{m}{1},\ldots,\id{m}{j-1} } \lfcP{ n+1-j }{ \id{m}{n},\ldots,\id{m}{j} }
	}
Let $t$ be a variable.
By substituting $(m_1-xt,m_n-yt)$ for $(m_1,m_n)$ and a modification, we obtain
	\envMLinePd[2_Lem1P_Eq1]
	{
		\lfcP{n}{ m_1-xt, m_2,\ldots,m_{n-1},m_n-yt } + \mo^{n} \lfcP{n}{ m_n-yt, m_{n-1},\ldots,m_2,m_1-xt }
		\\
	}
	{
		\tpSm{j_1,j_2\geq1}{j_1+j_2=n} \mo^{j_2-1} \lfcP{j_1}{ m_1-xt,m_2,\ldots,m_{j_1} } \lfcP{j_2}{ m_n-yt,m_{n-1}\ldots,m_{j_1+1} }
	}
We set first and second terms of left hand side with $A_1$ and $A_2$ respectively, 
	and the right hand side with $B$.
For each $A_1,A_2$ and $B$ up to $z^{m_1+\cdots+m_n}$,
	we will consider taking the sum of positive integers $m_1,\ldots,m_n\geq1$ and differentiating at $t=0$ with $l-n$ times.
	\envLineSePd
	{
		\bpDer[a]{t}_{t=0}^{l-n} \Sm{m_1,\ldots,m_n\geq1} \pw{z}{m_1+\cdots+m_n} A_1
	}
	{
		 \bpDer[a]{t}_{t=0}^{l-n} \Sm{m_1,\ldots,m_n\geq1} \pw{z}{m_1+\cdots+m_n} \lfcP{n}{ \id{m}{1}-x t,\id{m}{2},\ldots,\id{m}{n-1},\id{m}{n}-y t }
	}
	{
		\bpDer[a]{t}_{t=0}^{l-n} \Sm{m_1,\ldots,m_n\geq1}  \opF{\pw{z}{m_1+\cdots+m_n}}{m_1+\cdots+m_n - (x+y)t }  \PdT{i=2}{n} \opF{1}{ m_i+\cdots+m_n - yt }
	}
	{
		\bpDer[a]{t}_{t=0}^{l-n} \Sm{m_1>\ldots>m_n>0} \opF{ \pw{z}{m_1} }{m_1 - (x+y)t }  \PdT{i=2}{n} \opF{1}{ m_i - yt }
	}
	{
		(l-n)! \Sm{m_1>\ldots>m_n>0} \tpSm{ l_1, l_2,\ldots,l_n\geq0 }{l_1+\cdots+l_n=l-n} \pwR{x+y}{l_1} \pw{y}{l_2+\cdots+l_n}
		\opF{ \pw{z}{m_1} }{ \pw{m_1}{l_1+1}\cdots\pw{m_1}{l_n+1} }
	}
	{
		(l-n)! \tpSm{ l_1, l_2,\ldots,l_n\geq1 }{l_1+\cdots+l_n=l}  \Sm{m_1>\ldots>m_n>0} \pwR{x+y}{l_1-1} \pw{y}{l-l_1-(n-1)} 
		\opF{ \pw{z}{m_1} }{ \pw{m_1}{l_1}\cdots\pw{m_1}{l_n} }
	}
	{
		(l-n)! \lsmMP{ l }{ n }{ x+y,y }{ z }
	}
We can calculate those for $A_2$ and $B$ similarly, and we will get 
	\envHLineCm
	{
		\mo^{n} \bpDer[a]{t}_{t=0}^{l-n} \Sm{m_1,\ldots,m_n\geq1} \pw{z}{m_1+\cdots+m_n} A_2
	}
	{
		\mo^{n} (l-n)! \lsmMP{ l }{ n }{ y+x,x }{ z }
	}
	and
	\envLinePd
	{
		\bpDer[a]{t}_{t=0}^{l-n} \Sm{m_1,\ldots,m_n\geq1} \pw{z}{m_1+\cdots+m_n} B
	}
	{
		(l-n)! \tpSm{j_1,j_2\geq1}{j_1+j_2=n} \mo^{j_2-1} \tpSm{k_1,k_2\geq0}{k_1+k_2=l-n}
		\pw{x}{k_1} \pw{y}{k_2} \lgP{k_1+1,\pwB{1}{j_1-1}}{z} \lgP{k_2+1,\pwB{1}{j_2-1}}{z}
	}
These with \refEq{2_Lem1P_Eq1} prove the lemma.
}

In order to derive identities among multiple zeta values from Lemma \ref{2_Lem1},
	we need asymptotic properties of multiple polylogarithms $\lgP{l_1,\ldots,l_n}{z}$ with $l_1=1$ as $z\nearrow1$.
We do not have to consider the case of $l_1\geq2$ because 
	\envM{\Lim{z\nearrow1} \lgP{l_1,\ldots,l_n}{z} }{ \fcZ{l_1,\ldots,l_n} } if $l_1\geq2$.
\begin{lemma}[cf. \text{\cite[Section 2]{IKZ06}}]\label{2_Lem2}
Let $n$ be a positive integer.
For any multi-index $(l_1,\ldots,l_n)$ of positive integers with $l_1=1$,
	there is a polynomial $\idFc{f}{l_1,\ldots,l_n}{t}$ with real coefficients and a positive real number $J>0$ such that
	\envHLineCm
	{
		\lgP{l_1,\ldots,l_n}{z}
	}
	{
		\idFc[b]{f}{l_1,\ldots,l_n}{-\lgg[n]{1-z}} + \sbLandau[a]{(1-z) \pwR[b]{\lgg[n]{1-z}}{J}}		\qquad	(z\nearrow1)
	}
	where $\sbLandau[]{}$ denotes the Landau symbol.
\end{lemma}
\begin{remark}\label{2_Rem1}
In \cite{IKZ06},
	every polynomial $\idFc{f}{l_1,\ldots,l_n}{t}$ is given more concretely.
However we need only the fact that it is a polynomial with real coefficients in this paper.
\end{remark}
By Lemma \ref{2_Lem1} and \ref{2_Lem2}, we obtain the formula below,
	which with $x=y=1$ gives the formula of Eie, Yang and Ong \cite[Proposition 8]{EYO10}.
\begin{proposition}[cf. \text{\cite[Proposition 8]{EYO10}}] \label{2_Prop1}
For integers $l,n$ with $l>n\geq1$ and complex numbers $x,y$, we have
	\envLineCm[2_Prop1_Eq1]
	{
		\lsmZV{l}{n}{x+y,y} +\mo^{n}\lsmZV{l}{n}{y+x,x} - \bkR{ \mo^{n}\pw{x}{l-n}+\pw{y}{l-n} }\fcZ{l}
	}
	{
		\tpSm{j_1,j_2\geq1}{j_1+j_2=n} \mo^{j_2-1} \tpSm{k_1,k_2\geq1}{k_1+k_2=l-n}
		\pw{x}{k_1} \pw{y}{k_2} \fcZ{k_1+1,\pwB{1}{j_1-1}} \fcZ{k_2+1,\pwB{1}{j_2-1}}
		\lnAHs{150}
		+ 
		\bkR{ \pw{x}{l-n}+\mo^{n}\pw{y}{l-n} }\fcZ{l-n+1,\pwB{1}{n-1}}
	}
	where empty sum, occurring when $n=1$ or $l-n=1$, means $0$.
\end{proposition}
\envProof{
In \refEq{2_Lem1_Eq1}, the sum of the divergent terms of the left hand side as $z\nearrow1$ is 
	\envPLineCm
	{
		\bkR{ \mo^{n}\pw{x}{l-n}+\pw{y}{l-n} } \tpSm{l_2,\ldots,l_n\geq1}{l_2+\cdots+l_n=l-1} \lgP{1,l_2,\ldots,l_n}{z}
	}
	and that of the right hand side is
	\envLinePd
	{
		\tpSm{j_1,j_2\geq1}{j_1+j_2=n} \mo^{j_2-1} 
		\bkR{ \pw{x}{l-n}\lgP{l-n+1,\pwB{1}{j_1-1}}{z} \lgP{1,\pwB{1}{j_2-1}}{z} + \pw{y}{l-n}\lgP{1,\pwB{1}{j_1-1}}{z} \lgP{l-n+1,\pwB{1}{j_2-1}}{z} } 
	}
	{
		\bkR{ \mo^{n}\pw{x}{l-n}+\pw{y}{l-n} } \tpSm{j_1,j_2\geq1}{j_1+j_2=n} \mo^{j_2-1} \lgP{1,\pwB{1}{j_1-1}}{z} \lgP{l-n+1,\pwB{1}{j_2-1}}{z}
	}
By virtue of Lemma \ref{2_Lem2}, there is a polynomial $\Fc{f}{t}=\SmT{j=0}{m} a_j t^j$ 
	with real coefficients and a positive real number $J>0$ such that
	\envLinePd
	{
		\lsmZV{l}{n}{x+y,y} +\mo^{n}\lsmZV{l}{n}{y+x,x} 
	}
	{
		\tpSm{j_1,j_2\geq1}{j_1+j_2=n} \mo^{j_2-1} \tpSm{k_1,k_2\geq1}{k_1+k_2=l-n}
		\pw{x}{k_1} \pw{y}{k_2} \fcZ{k_1+1,\pwB{1}{j_1-1}} \fcZ{k_2+1,\pwB{1}{j_2-1}}		
		\lnAHs{30}
		+ 
		\bkR{ \mo^{n}\pw{x}{l-n}+\pw{y}{l-n} } \Fc[b]{f}{-\lgg[n]{1-z}} + \sbLandau[a]{(1-z) \pwR[b]{\lgg[n]{1-z}}{J}}			\qquad	(z\nearrow1)
	}
Since the right hand side must be converged, we see that $a_i=0\ (i=1,\ldots,m)$ and
	\envMLinePd
	{
		\lsmZV{l}{n}{x+y,y} +\mo^{n}\lsmZV{l}{n}{y+x,x} 
		\\
	}
	{
		\tpSm{j_1,j_2\geq1}{j_1+j_2=n} \mo^{j_2-1} \tpSm{k_1,k_2\geq1}{k_1+k_2=l-n}
		\pw{x}{k_1} \pw{y}{k_2} \fcZ{k_1+1,\pwB{1}{j_1-1}} \fcZ{k_2+1,\pwB{1}{j_2-1}}	
		+ 
		\bkR{ \mo^{n}\pw{x}{l-n}+\pw{y}{l-n} } a_0			
	}
By substituting $(0,1)$ for $(x,y)$ and by the sum formula \refEq{1_Eq_SF}, we obtain 
	\envHLineCm
	{
		a_0
	}
	{
		\fcZ{l} + \mo^{n}\fcZ{l-n+1,\pwB{1}{n-1}}
	}
	which completes the proof.
}
\begin{remark}\label{2_Rem2}
For any integer $r$ with $0\leq r\leq l-n$, 
	we define a real number $\lcoeSmZV{l}{n}{r}$ by the $\setQ$-linear combination of multiple zeta values
	\envOTLineDefPd
	{
		\lcoeSmZV{l}{n}{r}
	}
	{
		\tpSm{l_1\geq 1+ r \atop l_2,\ldots,l_n\geq1}{l_1+\cdots+l_n=l} \binom{l_1-1}{r} \fcZ{l_1,\ldots,l_n}
	}
	{
		\tpSm{l_1,\ldots,l_n\geq1}{l_1+\cdots+l_n=l-r} \binom{l_1+r-1}{r} \fcZ{l_1+r,l_2,\ldots,l_n}
	}
Here $\fcZ{1,l_2,\ldots,l_n}$ stands for $0$, thus $\lcoeSmZV{l}{n}{0} = \fcZ{l}$ by the sum formula \refEq{1_Eq_SF}.
Since
	\envHLine
	{
		\smZV{l}{n}{x+y,y}
	}
	{
		\SmT{j=0}{l-n} \pwR{x+y}{j}\pw{y}{l-n-j} \tpSm{l_2,\ldots,l_n\geq1}{l_2+\cdots+l_n=l-j-1} \fcZ{j+1,l_2,\ldots,l_n}
	}
	and
	\envOTLineCm
	{
		\SmT{j=0}{l-n} \pwR{x+y}{j}\pw{y}{l-n-j}
	}
	{
		\SmT{j=0}{l-n} \pw{y}{l-n-j} \SmT{i=0}{j}\binom{j}{i} \pw{x}{i} \pw{y}{j-i}
	}
	{
		\SmT{i=0}{l-n} \pw{x}{i} \pw{y}{l-n-i} \SmT{j=i}{l-n} \binom{j}{i} 
	}
	we have
	\envHLinePd
	{
		\lsmZV{l}{n}{x+y,y}
	}
	{
		\SmT{r=0}{l-n} \pw{x}{r} \pw{y}{l-n-r} \lcoeSmZV{l}{n}{r}
	}
By comparing of each coefficient of $x^{r}y^{l-n-r}\ (r=1,\ldots,l-n-1)$ in \refEq{2_Prop1_Eq1}, we get
	\envHLineCm
	{
		\lcoeSmZV{l}{n}{r} + \mo^{n} \lcoeSmZV{l}{n}{l-n-r}
	}
	{
		\tpSm{j_1,j_2\geq1}{j_1+j_2=n} \mo^{j_2-1}  \fcZ{r+1,\pwB{1}{j_1-1}} \fcZ{l-n-r+1,\pwB{1}{j_2-1}}
	}
	which is the same as the formula of Arakawa and Kaneko \cite[Corollary 11]{AK99} with $(r,m,k)=(r,l-n-r,n)$
	by virtue of the duality formula $\fcZ{p+1,\pwB{1}{q-1}} = \fcZ{q+1,\pwB{1}{p-1}}\ (p,q\geq1)$.
Proposition \ref{2_Prop1} and \cite[Corollary 11]{AK99} can be proved by each other, or are equivalent,
	since the remaining coefficients of $x^{l-n}$ and $y^{l-n}$ in  \refEq{2_Prop1_Eq1} mean the sum formula \refEq{1_Eq_SF}.
(We note that the order of the multi-index of the multiple zeta values in this paper and \cite{AK99} is reverse.)
\end{remark}

We are in a position to prove Theorem \ref{1_Thm1} now.

\envProof[Theorem \ref{1_Thm1}]{
We see from \cite{Aomoto90,BBB97,Drinfeld90} that
	\envHLinePd
	{
		\opF{ \fcGam{1-X}\fcGam{1-Y} }{ \fcGam{1-X-Y} }
	}
	{
		1 - \Sm{k,j\geq1} \fcZ{k+1,\pwB{1}{j-1}} \pw{X}{k} \pw{Y}{j}
	}
Thus
	\envLineThCm
	{
		\opF{ \fcGam{1-xX}\fcGam{1-Y} }{ \fcGam{1-xX-Y}} \opF{ \fcGam{1-yX}\fcGam{1+Y} }{ \fcGam{1-yX+Y} }
	}
	{
		1- \Sm{k,j\geq1} \pw{x}{k} \fcZ{k+1,\pwB{1}{j-1}} \pw{X}{k} \pw{Y}{j} - \Sm{k,j\geq1} \mo^{j} \pw{y}{k} \fcZ{k+1,\pwB{1}{j-1}} \pw{X}{k} \pw{Y}{j}
		\lnAHs{10}
		+
		\bkR[g]{ \Sm{k_1,j_1\geq1} \pw{x}{k_1} \fcZ{k_1+1,\pwB{1}{j_1-1}} \pw{X}{k_1} \pw{Y}{j_1} }
		\bkR[g]{ \Sm{k_2,j_2\geq1} \mo^{j_2} \pw{y}{k_2} \fcZ{k_2+1,\pwB{1}{j_2-1}} \pw{X}{k_2} \pw{Y}{j_2} }
	}
	{
		1
		-
		\Sm{l,n\geq1} \pw{X}{l} \pw{Y}{n}
		\bkS[G]{
			\bkR{ \pw{x}{l} + \mo^{n}\pw{y}{l} } \fcZ{l+1,\pwB{1}{n-1}}
			\lnAHs{80}
			+
			\tpSm{j_1,j_2\geq1}{j_1+j_2=n} \tpSm{k_1,k_2\geq1}{k_1+k_2=l}	
			\mo^{j_2-1} \pw{x}{k_1} \pw{y}{k_2} \fcZ{k_1+1,\pwB{1}{j_1-1}} \fcZ{k_2+1,\pwB{1}{j_2-1}}
		}
	}
	which together with Proposition \ref{2_Prop1} gives the first equation in Theorem \ref{1_Thm1}.
The second follows from the first and the well-known Taylor expansion
	\envMCm
	{
		\lgg{\fcGam{1-z}}
	}
	{
		\gamma z + \SmT[l]{m=2}{\infty} \opF[d]{\fcZ{m}}{m}z^m 
	}
	where $\gamma$ is the Euler constant.
}

\section{\sectThree} \label{sectThree}
In order to prove Theorem \ref{1_Thm2}, 
	we first give a different expression from the formula in Theorem \ref{1_Thm1} by use of the function 
	\envHLineDefPd
	{
		\lgfcZV{x,y,z}
	}
	{
		\Sm{n\geq1} \pw{z}{n} \Sm{m_1>\cdots>m_n>0}  \opF{ x }{ m_1 \bkR{m_1-x}\bkR{m_2-y}\cdots\bkR{m_n-y} }
	}
Since
	\envHLineFCm
	{
		\lgfcZV{xX,yX,Y}
	}
	{
		\Sm{n\geq1} \pw{Y}{n} \Sm{m_1>\cdots>m_n>0} \opF{x}{m_1}
		\opF{ X }{ \bkR{m_1-xX}\bkR{m_2-yX}\cdots\bkR{m_n-yX} }
	}
	{
		\Sm{n\geq1} \pw{Y}{n} \Sm{m_1>\cdots>m_n>0} \opF{x}{m_1} 
		\Sm{l\geq0} \pw{X}{l+1} \tpSm{l_1,\ldots,l_n\geq0}{l_1+\cdots+l_n=l}
		\opF{ \pw{x}{l_1} \pw{y}{l_2+\cdots+l_n} }{ \pw{m_1}{l_1+1}\pw{m_2}{l_2+1}\cdots\pw{m_n}{l_n+1} }
	}
	{
		\Sm{l,n\geq1} \pw{X}{l} \pw{Y}{n}  \tpSm{l_1\geq2,l_2,\ldots,l_n\geq1}{l_1+\cdots+l_n=l+n} 
		\pw{x}{l_1-1} \pw{y}{l_2+\cdots+l_n-(n-1)} \fcZ{l_1,\ldots,l_N}
	}
we see that $\lgfcZV{xX,yX,Y}$ is the generating function of $\lsmZV{l}{n}{x,y}$'s, or
	\envOTLinePd[3_Plane_Eq1]
	{
		\lgfcZV{xX,yX,Y}
	}
	{
		\Sm{l,n\geq1} \lsmZV{l+n}{n}{x,y} \pw{X}{l} \pw{Y}{n}
	}
	{
		\Sm{l>n\geq1} \lsmZV{l}{n}{x,y} \pw{X}{l-n} \pw{Y}{n}
	}
Therefore, by \refEq{3_Plane_Eq1} and \refEq{1_Thm1_Eq1} with $(xX,yX,Y) = (x,y,Y)$, 
	we obtain a functional equation for $\lgfcZV{x,y,Y}$ and $\fcGam{z}$, or the different expression
	\envLinePd[3_Plane_Eq2]
	{
		\lgfcZV{x+y,y,Y} + \lgfcZV{y+x,x,-Y} - \Sm{l,n\geq1} \fcZ{l+n} \bkR{ \pw{x}{l} \pwR{-Y}{n} + \pw{y}{l} \pw{Y}{n} }
	}
	{
		1 - \opF{ \fcGam{1-x}\fcGam{1-Y} \fcGam{1-y}\fcGam{1+Y} }{ \fcGam{1-x-Y}\fcGam{1-y+Y} }
	}
This implies the following formula.
\begin{proposition}\label{3_Prop1}
We have
	\envLinePd[3_Prop1_Eq1]
	{
		\lgfcZV{X+Y,Y,Y-X} + \lgfcZV{Y+X,X,X-Y} 
	}
	{
		1 - \opF{ \pi (Y-X) }{ \tfcS{\pi (Y-X)} } + \Sm{l,n\geq1} \fcZ{l+n} \bkR{ \pw{X}{l} \pwR{X-Y}{n}  + \pw{Y}{l} \pwR{Y-X}{n}  }
	}
\end{proposition}
\envProof{
By substituting $(X,X+Y)$ for $(x,y)$ in \refEq{3_Plane_Eq2}, we get
	\envLineThPd
	{
		\lgfcZV{2X+Y,X+Y,Y} + \lgfcZV{2X+Y,X-Y} - \Sm{l,n\geq1} \fcZ{l+n} \bkR{ \pw{X}{l} \pwR{-Y}{n} + \pwR{X+Y}{l} \pw{Y}{n} }
	}
	{
		1 - \fcGam{1-Y}\fcGam{1+Y}
	}
	{
		1 - \opF{ \pi Y }{ \tfcS{\pi Y } }
	}
We complete the proof of the proposition by replacing $Y$ with $Y-X$.
}
We prove Theorem \ref{1_Thm2}.
\envProof[Theorem \ref{1_Thm2}]{
Since
	\envLineFiCm
	{
		\Sm{l,n\geq1} \fcZ{l+n} \bkR{ \pw{X}{l} \pwR{X-Y}{n}  + \pw{Y}{l} \pwR{Y-X}{n} }
	}
	{
		\SmT{m=2}{\infty} \fcZ{m} \tpSm{l,n\geq1}{l+n=m} \bkR{ \pw{X}{l} \pwR{X-Y}{n}  + \pw{Y}{l} \pwR{Y-X}{n} }
	}
	{
		\SmT{m=2}{\infty} \fcZ{m} \bkR{ X(X-Y)\opF{ X^{m-1} - (X-Y)^{m-1} }{Y} + Y(Y-X)\opF{ Y^{m-1} - (Y-X)^{m-1} }{X} }
	}
	{
		\SmT{m=2}{\infty} \fcZ{m} \bkS{ (Y-X)\bkR{ \opF{Y^m}{X}-\opF{X^m}{Y} } - (Y-X)^m\bkR{ \mo^m\opF{X}{Y} + \opF{Y}{X} } }
	}
	{
		\SmT{m=2}{\infty} \fcZ{m} \bkS{ (Y-X)^2 \opF{Y^{m+1}-X^{m+1}}{YX(Y-X)}  - (Y-X)^m\bkR{ \mo^m\opF{X}{Y} + \opF{Y}{X} } }
	}
	and
	\envOTLineCm
	{
		\opF{ \pi z }{ \tfcS{\pi z } }
	}
	{
		\opF{\pi z}{2} \bkR{ \tfcCT{\opF{z}{2}}+\tfcT{\opF{z}{2}} }
	}
	{
		1 + 2 \tpSmT{l=2}{l\ even}{\infty} \opF{2^{l-1}-1}{2^{l-1}} \fcZ{m} z^l
	}
	we obtain the following equations;
	\envMLineCm
	{
		\Sm{l,n\geq1} t^{l+n} \fcZ{l+n} \bkR{ \pw{x}{l} \pwR{x-y}{n}  + \pw{y}{l} \pwR{y-x}{n}  }
		\\
	}
	{
		\SmT{l=2}{\infty} t^l (y-x)^2\bkS{ \opF{y^{l+1}-x^{l+1}}{yx(y-x)} - (y-x)^{l-2}\bkR{ \mo^l\opF{x}{y} + \opF{y}{x} } } \fcZ{l}
	}
\vPack
	\envHLinePd
	{
		1 - \opF{ \pi (y-x)t }{ \tfcS{\pi (y-x)t} }
	}
	{
		- \SmT{l=2}{\infty} t^l (1+\mo^l) \opF{2^{l-1}-1}{2^{l-1}} (y-x)^l \fcZ{l}
	}
By \refEq{3_Plane_Eq1}, we also get
	\envMLinePd
	{
		\lgfcZV{(x+y)t,yt,(y-x)t} + \lgfcZV{(x+y)t,xt,(x-y)t} 
		\\
	}
	{
		\Sm{l\geq2} \pw{t}{l} \SmT{n=1}{l-1} \pwR{y-x}{n} \bkR{ \lsmZV{l}{n}{x+y,y} + \mo^{n}\lsmZV{l}{n}{y+x,x} }
	}
Calculating the coefficients of $\pw{t}{l}$'s of the above equations  and Proposition \ref{3_Prop1} prove the theorem.
}

\section{\sectFour} \label{sectFour}
In this section, we give some identities among multiple zeta values with small depths between $2$ and $4$ by use of the theorems.
First we obtain them by direct calculations of the coefficients of $Y^n\ (n=2,3,4)$ of \refEq{1_Thm1_Eq1} in Theorem \ref{1_Thm1}.

\begin{proposition}\label{4_Prop1}
Let $l$ be an integer, and $\delta_{i,j}$ Kronecker's delta which is $1$ if $i=j$ and $0$ otherwise.
We assume that empty sum is equal to $0$.
\mbox{}\\{\bf (i)} 
If $l\geq3$, 
	\envHLinePd
	{
		\lsmZV{l}{2}{x+y,y} + \lsmZV{l}{2}{y+x,x}
	}
	{
		\opF{l+1}{2}\bkR[b]{ \pw{x}{l-2} + \pw{y}{l-2} } \fcZ{l} 
		-
		\opF{1}{2} \tpSm{j_1,j_2\geq2}{j_1+j_2=l} \PdT{i=1}{2}\bkR[B]{ \pw{x}{j_i-1} - \pw{y}{j_i-1} } \fcZ{j_i}  
	}
{\bf (ii)}
If $l\geq4$, 
	\envHLinePd
	{
		\lsmZV{l}{3}{x+y,y} - \lsmZV{l}{3}{y+x,x}
	}
	{
		\opF{(l+1)(l-4)}{6}\bkR[b]{ \pw{x}{l-3} - \pw{y}{l-3} } \fcZ{l} 
		\lnAHs{0}
		-
		\opF{1}{2}\tpSm{j_1\geq2,j_2\geq3}{j_1+j_2=l} (j_2-1) \PdT{i=1}{2 } \bkR[B]{ \pw{x}{j_i-1-\delta_{i,2}} - \mo^{\delta_{i,2}} \pw{y}{j_i-1-\delta_{i,2}} } \fcZ{j_i} 
		\lnAHs{0}
		+
		\opF{1}{6} \tpSm{j_1,j_2,j_3\geq2}{j_1+j_2+j_3=l}  \PdT{i=1}{3}\bkR[B]{ \pw{x}{j_i-1} - \pw{y}{j_i-1} } \fcZ{j_i}
	}
{\bf (iii)}
If $l\geq5$, 
	\envHLinePd
	{
		\lsmZV{l}{4}{x+y,y} + \lsmZV{l}{4}{y+x,x}
	}
	{
		\opF{(l+1)(l^2-7l+18)}{24} \bkR[b]{ \pw{x}{l-4} + \pw{y}{l-4} } \fcZ{l}
		\lnAHs{0}
		-
		\opF{1}{6} \tpSm{j_1\geq2,j_2\geq4}{j_1+j_2=l} (j_2-1)(j_2-2)  \PdT{i=1}{2 } \bkR[B]{ \pw{x}{j_i-1-2\delta_{i,2}} - \pw{y}{j_i-1-2\delta_{i,2}} } \fcZ{j_i} 
		\lnAHs{0}
		-
		\opF{1}{8} \tpSm{j_1,j_2\geq3}{j_1+j_2=l}  \PdT{i=1}{2 } (j_i-1) \bkR[B]{ \pw{x}{j_i-2} + \pw{y}{j_i-2} } \fcZ{j_i} 
		\lnAHs{0}
		+
		\opF{1}{4} \tpSm{j_1,j_2\geq2,j_3\geq3}{j_1+j_2+j_3=l}  (j_3-1)
		\PdT{i=1}{3}\bkR[B]{ \pw{x}{j_i-1-\delta_{i,3}} - \mo^{\delta_{i,3}} \pw{y}{j_i-1-\delta_{i,3}} } \fcZ{j_i}  
		\lnAHs{0}
		-
		\opF{1}{24}  \tpSm{j_1,j_2,j_3,j_4\geq2}{j_1+j_2+j_3+j_4=l} \PdT{i=1}{4}\bkR[B]{ \pw{x}{j_i-1} - \pw{y}{j_i-1} } \fcZ{j_i}  
	}
\end{proposition}
\envProof{
We set the argument of the exponential function in \refEq{1_Thm1} with 
	\envHLineDefPd
	{
		\lfcF{x,y}{X,Y}
	}
	{
		\Sm{m\geq2} \fcZ{m}\opF{ \bkR{ \pw{x}{m}+\pw{y}{m} }\pw{X}{m}+(1+\mo^{m})\pw{Y}{m}-\pw[n]{x X+Y}{m} - \pw[n]{y X-Y}{m} }{m}
	}
Since
	\envLineFCm
	{
		\Sm{m\geq2} \fcZ{m}\opF{ \pw[n]{x X+Y}{m} + \pw[n]{y X-Y}{m} }{m}
	}
	{
		\Sm{m\geq2} \opF{\fcZ{m}}{m}
		\bkS{ \SmT{n=0}{m}\binom{m}{n}\pw{Y}{n}\pw[n]{xX}{m-n} + \SmT{n=0}{m}\binom{m}{n}\pw[n]{-Y}{n}\pw[n]{yX}{m-n} }
	}
	{
		\Sm{n\geq0}\pw{Y}{n} \Sm{m\geq\max\bkB{2,n}} \binom{m}{n}\opF{\fcZ{m}}{m} \bkR[B]{ \pw{x}{m-n} + \mo^{n}\pw{y}{m-n} } \pw{X}{m-n}
	}
	{
		\Sm{m\geq2} \opF{\fcZ{m}}{m} \bkR[B]{ \pw{x}{m} + \pw{y}{m} } \pw{X}{m}
		+
		\pw{Y}{1} \Sm{m\geq1} (m+1)\opF{\fcZ{m+1}}{m+1} \bkR[B]{ \pw{x}{m} - \pw{y}{m} } \pw{X}{m}
		\lnAHs{30}
		+
		\Sm{n\geq2}\pw{Y}{n} \bkS[G]{ \opF{\fcZ{n}}{n} \bkR[B]{ 1 + \mo^{n} } + \Sm{m\geq1} \binom{m+n}{n}\opF{\fcZ{m+n}}{m+n} \bkR[B]{ \pw{x}{m} + \mo^{n}\pw{y}{m} } \pw{X}{m} }
	}
we get
	\envHLinePd
	{
		\lfcF{x,y}{X,Y}
	}
	{
		- \Sm{n\geq1}\pw{Y}{n} \Sm{m\geq1} \pw{X}{m}  \binom{m+n}{n}\opF{\fcZ{m+n}}{m+n} \bkR[B]{ \pw{x}{m} + \mo^{n}\pw{y}{m} } 
	}
We will calculate the Taylor polynomials of $- \pwR{ \lfcF{x,y}{X,Y}}{j}\ (j=1,2,3,4)$ with degree $4$ at $Y=0$.
	\envHLineCEPdPt{\sim}
	{
		- \lfcF{x,y}{X,Y}
	}
	{
		Y \Sm{l\geq2} \pw{X}{l-1} \bkR[B]{ \pw{x}{l-1} - \pw{y}{l-1} } \fcZ{l}  
		\lnAHs{0}
		+
		\opF{Y^2}{2} \Sm{l\geq3} \pw{X}{l-2}  (l-1) \bkR[B]{ \pw{x}{l-2} + \pw{y}{l-2} } \fcZ{l} 
		\lnAHs{0}
		+
		\opF{Y^3}{6} \Sm{l\geq4} \pw{X}{l-3}  (l-1)(l-2) \bkR[B]{ \pw{x}{l-3} - \pw{y}{l-3} } \fcZ{l} 
		\lnAHs{0}
		+
		\opF{Y^4}{24} \Sm{l\geq5} \pw{X}{l-4} (l-1)(l-2)(l-3) \bkR[B]{ \pw{x}{l-4} + \pw{y}{l-4} } \fcZ{l}  
	}
	{
		- \pwR{\lfcF{x,y}{X,Y} }{2}
	}
	{
		-
		Y^2 \Sm{l\geq4} \pw{X}{l-2} \tpSm{j_1,j_2\geq2}{j_1+j_2=l} \PdT{i=1}{2}\bkR[B]{ \pw{x}{j_i-1} - \pw{y}{j_i-1} } \fcZ{j_i}  
		\lnAHs{0}
		-
		Y^3 \Sm{l\geq5} \pw{X}{l-3} 
		\tpSm{j_1\geq2,j_2\geq3}{j_1+j_2=l} (j_2-1) \PdT{i=1}{2 } \bkR[B]{ \pw{x}{j_i-1-\delta_{i,2}} - \mo^{\delta_{i,2}} \pw{y}{j_i-1-\delta_{i,2}} } \fcZ{j_i} 
		\lnAHs{0}
		-
		Y^4\bkB[G]{
			\opF{1}{3} \Sm{l\geq6} \pw{X}{l-4} \tpSm{j_1\geq2,j_2\geq4}{j_1+j_2=l} (j_2-1)(j_2-2)  
			\PdT{i=1}{2 } \bkR[B]{ \pw{x}{j_i-1-2\delta_{i,2}} - \pw{y}{j_i-1-2\delta_{i,2}} } \fcZ{j_i} 
			\lnAHs{40}
			+
			\opF{1}{4} \Sm{l\geq6} \pw{X}{l-4} \tpSm{j_1,j_2\geq3}{j_1+j_2=l}  \PdT{i=1}{2 } (j_i-1) \bkR[B]{ \pw{x}{j_i-2} + \pw{y}{j_i-2} } \fcZ{j_i} 
		}
	}
	{
		-\pwR{ \lfcF{x,y}{X,Y} }{3}
	}
	{
		Y^3 \Sm{l\geq6} \pw{X}{l-3} \tpSm{j_1,j_2,j_3\geq2}{j_1+j_2+j_3=l}  \PdT{i=1}{3}\bkR[B]{ \pw{x}{j_i-1} - \pw{y}{j_i-1} } \fcZ{j_i}  
		\lnAHs{0}
		+
		Y^4 \opF{3}{2} \Sm{l\geq7} \pw{X}{l-4} \tpSm{j_1,j_2\geq2,j_3\geq3}{j_1+j_2+j_3=l}  (j_3-1)
		\PdT{i=1}{3}\bkR[B]{ \pw{x}{j_i-1-\delta_{i,3}} - \mo^{\delta_{i,3}} \pw{y}{j_i-1-\delta_{i,3}} } \fcZ{j_i}  
	}
	{
		-\pwR{ \lfcF{x,y}{X,Y} }{4}
	}
	{
		-
		Y^4 \Sm{l\geq8} \pw{X}{l-4} \tpSm{j_1,j_2,j_3,j_4\geq2}{j_1+j_2+j_3+j_4=l} \PdT{i=1}{4}\bkR[B]{ \pw{x}{j_i-1} - \pw{y}{j_i-1} } \fcZ{j_i}  
	}
On the other hand, we see from Theorem \ref{1_Thm1} that
	\envLineThPd
	{
		\Sm{l>n\geq1} \bkS[g]{ \lsmZV{l}{n}{x+y,y} + \mo^{n} \lsmZV{l}{n}{y+x,x} - \bkR[b]{ \mo^{n}\pw{x}{l-n} + \pw{y}{l-n} } \fcZ{l} } \pw{X}{l-n} \pw{Y}{n}
	}
	{
		1 - \exx[b]{ \lfcF{x,y}{X,Y} }
	}
	{
		- \pw{ \lfcF{x,y}{X,Y} }{1} - \opF{ \pw{ \lfcF{x,y}{X,Y} }{2} }{ 2 } - \opF{ \pw{ \lfcF{x,y}{X,Y} }{3} }{ 6 } - \opF{ \pw{ \lfcF{x,y}{X,Y} }{4} }{ 24 }  + \sbLandau{Y^5}
	}
We obtain the required equations by compering the coefficients of $X^lY^n$'s.
}
\begin{remark}\label{4_Rem1}
We can not derive the formula \refEq{1_Eq_WSFp} of Gangl, Kaneko and Zagier from (i) in Proposition \ref{4_Prop1} directly,
	but we can do it by use of the Euler sum formula \refEq{1_Eq_eulerSF} and the harmonic relations among the double zeta values
	\envHLine
	{
		\fcZ{m}\fcZ{n}
	}
	{
		\fcZ{m,n} + \fcZ{n,m} + \fcZ{m+n} 
	}
	where $m,n\geq2$.
On the other hand, we obtain the weighted Euler sum formula \refEq{1_Eq_eulerWSF} of Ohno and Zudilin 
	by (i) with $x=y=1$ in Proposition \ref{4_Prop1} directly.
\end{remark}
As applications of \refH{Theorem}{1_Thm2} and \refH{Proposition}{4_Prop1}, we will give identities including the weighted sums 
	\envPLine
	{
		\tpSm{l_1\geq2,\l_2,l_3,l_4\geq1}{l_1+l_2+l_3+l_4=l} \pw{2}{l_1-1} \fcZ{l_1,l_2,l_3,l_4}
	}
	which appear in the weighted sum formula of Guo and Xie \cite[Theorem 1.1]{GX09} in the case of depth $4$,
	\envHLinePd
	{
		l\fcZ{l}
	}
	{
		\tpSm{l_1\geq2,l_2,l_3,l_4\geq1}{l_1+l_2+l_3+l_4=l}
		\bkS[g]{ \pw{2}{l_1-1} + \bkR{\pw{2}{l_1-1}-1} \bkR{ \pw{2}{l_2-1} + \pw{2}{l_2+l_3-1} } } \fcZ{l_1,l_2,l_3,l_4}
	}
Unlike the formula of Guo and Xie, 
	our identities include either weighted sums of multiple zeta values with lower depths $2,3$
	or sums of two products of single zeta values.
	
\begin{proposition}\label{4_Prop2}
Let $l$ be an integer with $l\geq5$. 
\mbox{}\\{\bf (i)} 
It holds that
	\envHLinePd
	{
		\opF{ (l+1)(l^2-l+6) }{48} \fcZ{l}
	}
	{
		\tpSm{l_1\geq2,\l_2,l_3,l_4\geq1}{l_1+l_2+l_3+l_4=l} \pw{2}{l_1-1} \fcZ{l_1,l_2,l_3,l_4}
		\lnAHs{10}
		+
		\opF{1}{2} \tpSm{l_1\geq2,\l_2,l_3\geq1}{l_1+l_2+l_3=l} (l-l_1-2)\pw{2}{l_1-1} \fcZ{l_1,l_2,l_3}
		\lnAHs{10}
		+
		\opF{1}{8} \tpSm{l_1\geq2,\l_2,\geq1}{l_1+l_2=l} (l-l_1-1)(l-l_1-2) \pw{2}{l_1-1} \fcZ{l_1,l_2}
	}
\mbox{}\\{\bf (ii)} 
It holds that
	\envHLinePd
	{
		\opF{(l+1)(l^2-7l+18)}{24} \fcZ{l}
	}
	{
		\tpSm{l_1\geq2,\l_2,l_3,l_4\geq1}{l_1+l_2+l_3+l_4=l} \pw{2}{l_1-1} \fcZ{l_1,l_2,l_3,l_4}
		+
		\opF{1}{4} \tpSm{j_1,j_2\geq3}{j_1+j_2=l}  \PdT{i=1}{2 } (j_i-1) \fcZ{j_i} 
	}
\end{proposition}
We prepare a notation and a lemma to prove the proposition. 
Let $\lsmZVd{l}{n}{p,q}{x,y}$ be the derivative function $\displaystyle\pDerT{\ptl x^p \ptl y^q}{p+q} \lsmZV{l}{n}{x,y}$
	for any pairs $(p,q)$ of nonnegative integers and $(l,n)$ of positive integers with $l>n$.
By \refEq{1_Def_wsmMZV}, it is expressed as
	\envHLineCm[4_Def_wsmMZVd]
	{\hspace{-20pt}
		\lsmZVd{l}{n}{p,q}{x,y}
	}
	{
		\tpSm{l_1\geq2,  l_2,\ldots,l_n\geq1}{l_1+\cdots+l_n=l} 
		\sbPochD{l_1-1}{p}\sbPochD{l-l_1-(n-1)}{q} \pw{x}{l_1-1-p}\pw{y}{l-l_1-(n-1)-q} \fcZ{l_1,\ldots,l_n}
	}
	where $\sbPochD{x}{n}:=x(x-1)\cdots(x-n+1)$ denotes the descending Pochhammer symbol,
	in particular, $\sbPochD{x}{0}=1$ and $\sbPochD{m}{n} = 0$ for integers $m,n$ with $m<n$.
The lemma deals with special values of derivative functions.
\begin{lemma}\label{4_Lem1}
Let $l,n$ be integers satisfying $l\geq5, n\geq3$, and
	$\Psi$ the differential operator 
	\envHLineDefPd
	{
		\Psi
	}
	{
		\bkR{ \pDerT{\ptl y^2}{2} - \pDerT{\ptl x\ptl y}{2} }_{x=y=1}
	}
\mbox{}\\
{\bf (i)}
We have
	\envHLineCEPd
	{
		\Fc[g]{\Psi}{ (y-x)^n \bkR{ \lsmZV{l}{n}{x+y,y} + \mo^n \lsmZV{l}{n}{y+x,x} } }
	}
	{
		0
	}
	{
		\Fc[g]{\Psi}{ (y-x)^2 \bkR{ \lsmZV{l}{4}{x+y,y} + \lsmZV{l}{4}{y+x,x} } }
	}
	{
		8\lsmZV{l}{4}{2,1}
	}
	{
		\Fc[g]{\Psi}{ (y-x) \bkR{ \lsmZV{l}{3}{x+y,y} - \lsmZV{l}{3}{y+x,x} } }
	}
	{
		4\lsmZVd{l}{3}{0,1}{2,1}
	}
	{
		\Fc[g]{\Psi}{ \bkR{ \lsmZV{l}{2}{x+y,y} + \lsmZV{l}{2}{y+x,x} } }
	}
	{
		\lsmZVd{l}{2}{0,2}{2,1}
	}
{\bf (ii)}
We have
	\envHLineCFPd
	{
		\Fc[a]{\Psi}{ \opF{\pw{y}{l+1}-\pw{x}{l+1}}{yx(y-x)} } 
	}
	{
		\opF{ (l+1)(l^2-l+6) }{6}
	}
	{
		\Fc[G]{\Psi}{ \pwR{y-x}{l-2} \bkR[g]{ (1+\mo^l) \opF{2^{l-1}-1}{2^{l-1}} + \mo^{l}\opF{x}{y} + \opF{y}{x} } } 
	}
	{
		0
	}
	
\end{lemma}
\envProof{
We prove the equations in (i).
The first equation is evident because of $n\geq3$, 
	and the second is derived from $\Fc[b]{\Psi}{ (y-x)^2} = 4$.
The third and fourth follow since $\pDerT{\ptl y^2}{2} - \pDerT{\ptl x\ptl y}{2} = \bkR[b]{ \pDer{y} - \pDer{x} } \pDer{y}$ and
	\envMLineCm
	{
		\pDer{y} (y-x) \bkR{ \lsmZV{l}{3}{x+y,y} - \lsmZV{l}{3}{y+x,x} }
	}
	{
		\bkR{ \lsmZV{l}{3}{x+y,y} - \lsmZV{l}{3}{y+x,x} }
		\\
		+
		(y-x) \bkR{ \lsmZVd{l}{3}{1,0}{x+y,y} + \lsmZVd{l}{3}{0,1}{x+y,y} - \lsmZVd{l}{3}{1,0}{y+x,x} }
	}
\vPack[10]
	\envHLinePd
	{
		\pDer{y} \bkR{ \lsmZV{l}{2}{x+y,y} + \lsmZV{l}{2}{y+x,x} }
	}
	{
		\lsmZVd{l}{2}{1,0}{x+y,y} + \lsmZVd{l}{2}{0,1}{x+y,y} + \lsmZVd{l}{2}{1,0}{y+x,x}
	}

We verify the equations in (ii) next.
We see from 	
	\envM
	{
		\opF[s]{(y^{l-1}-x^{l-1})}{(y-x)}  
	}
	{
		\SmT[l]{j=0}{l-2} y^{j}x^{l-2-j}
	}
	that
	\envOTLinePd
	{
		\Fc[a]{\Psi}{ \opF{\pw{y}{l-1}-\pw{x}{l-1}}{y-x} }
	}
	{
		\SmT{j=1}{l-2} \bkR[b]{ j(j-1) -  j(l-2-j) }
	}
	{
		\SmT{j=1}{l-2} \bkR{ 2j^2 - (l-1)j }
	}
By use of the famous formula $\SmT[l]{j=1}{k}j = \opF[s]{k(k+1)}{2}$ and $\SmT[l]{j=1}{k}j^2 = \opF[s]{k(k+1)(2k+1)}{6}$, we get
	\envHLinePd
	{
		\Fc[a]{\Psi}{ \opF{\pw{y}{l-1}-\pw{x}{l-1}}{y-x} }
	}
	{
		\opF{(l-1)(l-2)(l-3)}{6}
	}
Thus we have
	\envHLineThCm
	{
		\Fc[a]{\Psi}{ \opF{\pw{y}{l+1}-\pw{x}{l+1}}{yx(y-x)} }
	}
	{
		\Fc[a]{\Psi}{ \opF{\pw{x}{l-1}}{y} + \opF{\pw{y}{l-1}}{x} + \opF{\pw{y}{l-1}-\pw{x}{l-1}}{y-x} } 
	}
	{
		(l+1) + (l-1)^2 + \opF{(l-1)(l-2)(l-3)}{6}
	\lnP{=}
		\opF{ (l+1)(l^2-l+6) }{6}
	}
	which proves the first equation in (ii).
The second is clear because of $l\geq5$.
}
\refH{Proposition}{4_Prop2} can be derived from \refH{Theorem}{1_Thm2}, \refH{Lemma}{4_Lem1} and \refH{Proposition}{4_Prop1} easily;
We see from \refH{Theorem}{1_Thm2} and \refH{Lemma}{4_Lem1} that
	\envHLineCm
	{
		8\lsmZV{l}{4}{2,1} + 4\lsmZVd{l}{3}{0,1}{2,1} + \lsmZVd{l}{2}{0,2}{2,1}
	}
	{
		\opF{ (l+1)(l^2-l+6) }{6} \fcZ{l}
	}
	which with \refEq{1_Def_wsmMZV} and \refEq{4_Def_wsmMZVd} verifies \refH{Proposition}{4_Prop2} (i).
The equation with $x=y=1$ in \refH{Proposition}{4_Prop1} (iii) proves \refH{Proposition}{4_Prop2} (ii).

\section*{Acknowledgements}
The author is grateful to Yasuo Ohno and Wadim Zudilin for their constructive comments,
	and to Masanobu Kaneko 
	for pointing out the relation between \refEq{2_Prop1_Eq1} and \cite[Corollary 11]{AK99}.


\begin{flushleft}
				\mbox{}\\ \qquad \mbox{}\\ \qquad
\majorName		\mbox{}\\ \qquad
\departmentName	\mbox{}\\ \qquad
\organizationName	\mbox{}\\ \qquad
\placeAddress		\mbox{}\\ \qquad
\emailAddress

\end{flushleft}

\end{document}